\documentclass[reqno]{amsart}
\usepackage{amssymb}
\usepackage{amsmath}
\usepackage{amsfonts}
\usepackage{amscd,amsthm}

\newtheorem{theorem}{Theorem}
\theoremstyle{plain}

\newtheorem{corollary}{Corollary}

\newtheorem{remark}{Remark}

\numberwithin{equation}{section}
\input{tcilatex}

\begin{document}

\begin{center}
\textbf{{\large Univalence criteria related with Ruscheweyh and S\u{a}l\u{a}%
gean derivatives and Loewner Chains }}

\bigskip

\textbf{\ Erhan Deniz}$^{\ast }$\textbf{\ and Halit Orhan}

\medskip

\textit{Department of Mathematics$,$ Faculty of Science}\\[0pt]
\textit{Atat\"{u}rk University$,$ TR-$25240$ Erzurum$,$ Turkey}\\[0pt]

\textbf{E-Mail: edeniz@atauni.edu.tr, \ horhan@atauni.edu.tr}\\[0pt]

\textbf{$^{\ast }$Corresponding Author} \\[0pt]

\bigskip

\textbf{Abstract}
\end{center}

\begin{quotation}
In this paper we obtain, by the method of Loewner chains, some sufficient
conditions for the analyticity and the univalence of the functions defined
by an integral operator. These conditions involves Ruscheweyh and S\u{a}l%
\u{a}gean derivative operator in the open unit disk. In particular cases, we
find the well-known conditions for univalency established by Becker \cite{Be}%
, Ahlfors \cite{Ah}, Kanas and Srivastava \cite{KaSr} and others for
analytic mappings $f:\mathcal{U}\rightarrow 
\mathbb{C}
.$ Also, we obtain the corresponding new, useful and simpler conditions for
this integral operator.
\end{quotation}



\section{\textbf{Introduction }}

Denote by $\mathcal{U}_{r}=\left\{ {z\in \mathbb{C}:\;\left\vert
z\right\vert <r}\right\} \quad (0<r\leqslant 1)$ the disk of radius $r$ and
let $\mathcal{U}=\mathcal{U}_{1}.$ Let $\mathcal{A}$ denote the class of
analytic functions in the open unit disk $\mathcal{U}$ which satisfy the
usual normalization condition:%
\begin{equation*}
f(0)={f}^{\prime }(0)-1=0.
\end{equation*}

By $\mathcal{S}$ the subclass of $\mathcal{A}$ consisting of functions $f(z)$
which are univalent in $\mathcal{U}$. These classes have been one of the
important subjects of research in Geometric Function Theory for a long time
(see \cite{SrOw} ). For the functions $f_{p}\;(p=1,2)$ given by%
\begin{equation*}
f_{p}(z)=z+\sum\limits_{k=2}^{\infty }{a_{k,p}z^{k}}\text{ \ \ }(p=1,2),
\end{equation*}

\noindent let $f_{1}\ast f_{2}$ denote the Hadamard product (or convolution)
of $f_{1}$ and $f_{2},$ defined by%
\begin{equation}
(f_{1}\ast f_{2})(z)=z+\sum\limits_{k=2}^{\infty }{a_{k,1}a_{k,2}z^{k}=}%
(f_{2}\ast f_{1})(z).  \label{eq1}
\end{equation}

Two of the most important and well-known univalence criteria for analytic
functions defined the open unit disk were obtained by Becker \cite{Be} and
Ahlfors \cite{Ah}. Becker and Ahlfors's works depends upon a clever use of
the Theory of Loewner chains and the generalized Loewner differential
equation. Extensions of these two criteria were given by Ruscheweyh \cite%
{Ru2}, Kanas and Srivastava \cite{KaSr}. Recently, Ovesea \cite{Ov}, Deniz
and Orhan \cite{DeOr} and Deniz \textit{et al}. \cite{DeRaor} obtained some
generalization of these univalence criterions. Furthermore, Pascu \cite{Pa}
and Raducanu \textit{et\ al}. \cite{RaOrDe} obtained some extensions of
Becker's univalence criterion.

In the present paper, we will study a number of new criteria for univalence
based upon the Ruscheweyh and S\u{a}l\u{a}gean derivative and for the
functions defined by the integral operator $\mathcal{F}_{\beta }(z).$ The
paper improve the work of Kanas and Srivastava \cite{KaSr} of extending
univalence criteria for analytic mappings. In special cases our univalence
conditions contain the results obtained by some of the authors are cited in
references. Our considerations are based on the Theory of Loewner chains.

\section{\textbf{Loewner chains and related theorem}}

Before proving our main theorem we need a brief summary of the method of
Loewner chains.

Let $\mathcal{L}(z,t)=a_{1}(t)z+a_{2}(t)z^{2}+...,$ $a_{1}(t)\neq 0,$ be a
function defined on $\mathcal{U}\times I$, where $I:=[0,\infty )$ and $%
a_{1}(t)$ is a complex-valued, locally absolutely continuous function on $I.$
$\mathcal{L}(z,t)$ is called a Loewner chain if $\mathcal{L}(z,t)$ satisfies
the following conditions;

\begin{enumerate}
\item[(i)] $\mathcal{L}(z,t)$ is analytic and univalent in $\mathcal{U}$ for
all $t\in I$

\item[(ii)] $\mathcal{L}(z,t)\prec \mathcal{L}(z,s)$ for all $0\leq t\leq
s<\infty $,
\end{enumerate}

where the symbol \textquotedblright $\prec $ \textquotedblright\ stands for
subordination. If $a_{1}(t)=e^{t}$ then we say that $\mathcal{L}(z,t)$ is a 
\textit{standard Loewner chain}.

In order to prove our main results we need the following theorem due to
Pommerenke \cite{Po1} (also see \cite{Po}). This theorem is often used to
find out univalency for an analytic function, apart from the theory of
Loewner chains;

\begin{theorem}
\label{T1}(see Pommerenke \cite{Po})\textbf{\ }\textit{Let }$\mathcal{L}%
(z,t)=a_{1}(t)z+a_{2}(t)z^{2}+...$\textit{\ be analytic in }$\mathcal{U}_{r}$%
\textit{\ for all }$t\in I.$ Suppose that;

\begin{enumerate}
\item[(i)] \textit{\ }$\mathcal{L}(z,t)$ is a \textit{locally absolutely
continuous function in\ the interval }$I,$\textit{\ and locally uniformly
with respect to }$\mathcal{U}_{r}.$\textit{\ }

\item[(ii)] $a_{1}(t)$ is a complex valued continuous function on $I$ such
that $a_{1}(t)\neq 0,$ $\left\vert {a_{1}(t)}\right\vert \rightarrow \infty $%
\textit{\ for }$t\rightarrow \infty $ and%
\begin{equation*}
\left\{ \frac{{\mathcal{L}(z,t)}}{a_{1}(t)}\right\} _{t\in I}
\end{equation*}%
\textit{\ forms a normal family of functions in }$\mathcal{U}_{r}.$

\item[(iii)] There exists an analytic function $p:\mathcal{U}\times
I\rightarrow 
\mathbb{C}
$ satisfying $\Re \left( {p(z,t)}\right) >0$\textit{\ for all }$z\in 
\mathcal{U},\;t\in I$ and 
\begin{equation}
z\frac{\partial \mathcal{L}(z,t)}{\partial z}=p(z,t)\frac{\partial \mathcal{L%
}(z,t)}{\partial t},\quad z\in \mathcal{U}_{r},\text{ }t\in I.  \label{1}
\end{equation}%
\textit{Then, for each }$t\in I,$\textit{\ the function }$\mathcal{L}(z,t)$%
\textit{\ has an analytic and univalent extension to the whole disk }$%
\mathcal{U}$ or the function $\mathcal{L}(z,t)$ is a Loewner chain.
\end{enumerate}
\end{theorem}

The equation (\ref{1}) is called the generalized Loewner differential
equation. Detailed information about Loewner chains and Pommerenke's theorem
can be found in Hotta \cite{Ho} and \cite{Ho1}.

\section{\textbf{Univalence criteria and the Ruscheweyh derivative}}

For a function $f\in \mathcal{A}$ the Ruscheweyh derivative operator \cite%
{Ru1}, $\mathcal{R}^{\lambda }:\mathcal{A}\rightarrow \mathcal{A}$ is
defined by%
\begin{equation}
\mathcal{R}^{\lambda }f(z)=\frac{z}{(1-z)^{\lambda +1}}\ast f(z)\text{ \ \ }%
\left( {\lambda >-1;\;z\in \mathcal{U}}\right) .  \label{eq2}
\end{equation}%
In particular, when $\lambda =n\;(n\in \mathbb{N}_{0}:=\mathbb{N}\cup
\{0\}), $ definition (\ref{eq2}) implies that%
\begin{equation}
\mathcal{R}^{n}f(z)=\frac{z}{n!}\frac{d^{n}}{dz^{n}}\{z^{n-1}f(z)\}.
\label{eq3}
\end{equation}%
It is obvious that%
\begin{eqnarray}
\mathcal{R}^{0}f(z) &=&f(z)  \label{eq4} \\
\mathcal{R}^{1}f(z) &=&z{f}^{\prime }(z)  \notag \\
\mathcal{R}^{2}f(z) &=&\frac{z}{2}\{2{f}^{\prime }(z)+z{f}^{\prime \prime
}(z)\},  \notag
\end{eqnarray}

\noindent and so on. Also we can write recurrence relationship as follows:%
\begin{equation*}
z[\mathcal{R}^{n}f(z){]}^{\prime }=(n+1)\mathcal{R}^{n+1}f(z)-n\mathcal{R}%
^{n}f(z).
\end{equation*}%
\qquad In this section, making use of the Theorem \ref{T1}, we obtain some
univalence criterions connected with the Ruscheweyh derivative operator for
an integral operator. The proofs are based on the theory of Loewner chains
(or Theorem \ref{T1}), essence of which is the construction of a suitable
Loewner chain.

\begin{theorem}
\label{T2} Let $f,g,h\in \mathcal{A}.$ Also \textit{let} $m\in \mathbb{R}%
_{+} $\textit{\ and }$\alpha ,$ $\beta ,$ $c$\textit{\ be complex numbers
such that}%
\begin{equation}
\mathit{\ }\alpha \neq 1,\text{ \ }c\neq -1,\;\ \left\vert \frac{1+c}{%
1-\alpha }-\frac{m+1}{2}\right\vert \leqslant \frac{m+1}{2},\text{ \ \ }%
\left\vert {\beta -\frac{m+1}{2}}\right\vert <\frac{m+1}{2}.  \label{eq411}
\end{equation}%
\textit{If the inequalities }%
\begin{equation}
\left\vert {\left( {\frac{(1+c){f}^{\prime }(z)}{[\mathcal{R}^{n}h(z){]}%
^{\prime }-\alpha }-1}\right) -\frac{m-1}{2}}\right\vert <\frac{m+1}{2}
\label{eq41}
\end{equation}%
\textit{and}%
\begin{equation*}
\left\vert {\left\vert z\right\vert ^{m+1}\left( {\frac{(1+c){f}^{\prime }(z)%
}{[\mathcal{R}^{n}h(z){]}^{\prime }-\alpha }-1}\right) }\right.
\end{equation*}%
\begin{equation}
\left. {+(1-\left\vert z\right\vert ^{m+1})\left[ {(\beta -1)\frac{z{g}%
^{\prime }(z)}{g(z)}+\frac{z[\mathcal{R}^{n}h(z){]}^{\prime \prime }}{[%
\mathcal{R}^{n}h(z){]}^{\prime }-\alpha }}\right] -\frac{m-1}{2}}\right\vert
\leqslant \frac{m+1}{2}  \label{eq5}
\end{equation}%
\textit{\ are true for all }$z\in \mathcal{U},$\textit{\ then the function }$%
\mathcal{F}_{\beta }(z)$ defined by%
\begin{equation}
\mathcal{F}_{\beta }(z)=\left[ {\beta \int\limits_{0}^{z}{g^{\beta -1}(u){f}%
^{\prime }(u)du}}\right] ^{1\diagup \beta }  \label{eq6}
\end{equation}%
\textit{is analytic and univalent in }$\mathcal{U},$\textit{\ where the
principal branch is intended.}
\end{theorem}

\begin{proof}
Let $a$ and $b$ be two positive real numbers such that $m=\frac{b}{a}.$ We
will prove that there exists a real number $r\in \left( {0,1}\right] $ such
that the function $\mathcal{L}:\mathcal{U}_{r}\times I\rightarrow \mathbb{C}%
, $ defined formally by%
\begin{eqnarray}
&&\mathcal{L}(z,t)  \label{eq7} \\
&=&\left\{ {\beta \int\limits_{0}^{e^{-at}z}{g^{\beta -1}(s){f}^{\prime
}(s)ds+}}\frac{\beta {(e^{bt}-e^{-at})}}{1+c}{zg^{\beta -1}(e^{-at}z)\left( {%
[\mathcal{R}^{n}h(e^{-at}z){]}^{\prime }-\alpha }\right) }\right\}
^{1\diagup \beta }  \notag
\end{eqnarray}%
\noindent is analytic in $\mathcal{U}_{r}$ for all $t\in I.$

Because $g\in \mathcal{A}$ the function%
\begin{equation*}
{\psi }(z)=\frac{g(z)}{z}
\end{equation*}%
\noindent is analytic in $\mathcal{U}$ and ${\psi }(0)=1$. Then there exist
a disk $\mathcal{U}_{r_{1}},$ $0<r_{1}\leqslant 1,$ in which ${\psi }(z)\neq
0$ for all $z\in \mathcal{U}_{r_{1}}.$ We denote by ${\psi }_{1}$ the
uniform branch of $\left( {\psi (z)}\right) ^{\beta -1}$ equal to 1 at
origin.

Consider the function%
\begin{equation*}
{\psi }_{2}(z,t)=\beta \int\limits_{0}^{e^{-at}z}{s^{\beta -1}\psi _{1}(s){f}%
^{\prime }(s)ds,}
\end{equation*}%
\noindent then we have%
\begin{equation*}
{\psi }_{2}(z,t)=z_{3}^{\beta }{\psi }(z,t)
\end{equation*}%
\noindent where ${\psi }_{3}$ is also analytic in $\mathcal{U}_{r_{1}}.$
Hence, the function%
\begin{equation*}
{\psi }_{4}(z,t)={\psi }_{3}(z,t)+\frac{\beta }{1+c}\left( {e^{bt}-e^{-at}}%
\right) e_{1}^{-a(\beta -1)t}{\psi }(e^{-at}z)\left( {[R{^{n}h}(e^{-at}z){]}%
^{\prime }-\alpha }\right)
\end{equation*}%
\noindent is analytic in $\mathcal{U}_{r_{1}}$ and%
\begin{equation*}
{\psi }_{4}(0,t)=e^{-a\beta t}\left[ {\frac{1+c-(1-\alpha )\beta }{1+c}+%
\frac{(1-\alpha )\beta }{1+c}(e^{(a+b)t})}\right] .
\end{equation*}%
We will prove that ${\psi }_{4}(0,t)\neq 0$ for all $t\in I.$ It is easy to
see that ${\psi }_{4}(0,0)=1.$ Suppose that there exists $t_{0}>0$ such that 
${\psi }_{4}(0,t_{0})=0.$ Then we obtain the equality $e^{(a+b)t_{0}}=\frac{%
1+c+(\alpha -1)\beta }{(\alpha -1)\beta }.$ This equality implies that $c<-1$
$($or $c>-1),$ which contradicts $\left\vert c\right\vert \leqslant 1.$ We
conclude that ${\psi }_{4}(0,t)\neq 0$ for all $t\in I.$ Therefore, there is
a disk $\mathcal{U}_{r_{2}},\;r_{2}\in \left( {0,r_{1}}\right] ,$ in which ${%
\psi }_{4}(z,t)\neq 0$ for all $t\in I.$ Then we can choose an uniform
branch of $\left[ {\psi _{4}(z)}\right] ^{1\diagup \beta }$ analytic in $%
\mathcal{U}_{r_{2}},$ denoted by ${\psi }_{5}(z,t).$

It follows from (\ref{eq7}) that the%
\begin{equation*}
\mathcal{L}(z,t)=z{\psi }_{5}(z,t)=a_{1}(t)z+a_{2}(t)z^{2}+...
\end{equation*}%
\noindent and thus, the function $\mathcal{L}(z,t)$ is analytic in $\mathcal{%
U}_{r_{2}}$.

We have%
\begin{equation*}
a_{1}(t)=e^{(\frac{-a\beta +a+b}{\beta })t}\left[ {\frac{1+c-(1-\alpha
)\beta }{1+c}e^{-(a+b)t}+}\frac{(1-\alpha )\beta }{1+c}\right] ^{1\diagup
\beta }.
\end{equation*}%
\noindent for which we consider the uniform branch equal to $1$ at the
origin. Since $\left\vert {\beta -\frac{m+1}{2}}\right\vert <\frac{m+1}{2}$
is equivalent with $\Re (\frac{1}{\beta })>\frac{1}{m+1}$ we have that%
\begin{equation*}
\underset{t\rightarrow \infty }{\lim }\left\vert {a_{1}(t)}\right\vert
=\infty .
\end{equation*}%
Moreover, $a_{1}(t)\neq 0$ for all $t\in I.$

After simple calculation we obtain, for each $z\in \mathcal{U}$%
\begin{eqnarray*}
&&\underset{t\rightarrow \infty }{\lim }\frac{\mathcal{L}(z,t)}{a_{1}(t)} \\
&=&\underset{t\rightarrow \infty }{\lim }z\left[ {\frac{1+c-(1-\alpha )\beta 
}{1+c}e^{-(a+b)t}+}\frac{(1-\alpha )\beta }{1+c}\right] ^{-1\diagup \beta }
\\
&&\times \left( 1+O({e^{-at}z})\right) ^{1\diagup \beta }\left\{ {e^{-(a+b)t}%
}+(1-{e^{-(a+b)t}})\frac{\beta }{1+c}\left[ 1-\alpha +O({e^{-at}z})\right]
\right\} ^{1\diagup \beta } \\
&=&z.
\end{eqnarray*}%
The limit function $\varphi (z)=z$ belongs to the family $\left\{ \mathcal{L}%
(z,t)\diagup a_{1}(t)\right\} ;$ then in every closed disk $\mathcal{U}%
_{r_{3}},$ $0<r_{3}<r_{2},$ there exists a constant $K=K(r_{3})$ such that%
\begin{equation*}
\left\vert {\frac{\mathcal{L}(z,t)}{a_{1}(t)}}\right\vert <K,\text{ \ \ }%
\forall z\in \mathcal{U}_{r_{3}},\;t\in I
\end{equation*}%
uniformly in this disk, provided that $t$ is sufficiently large. Then, by
Montel's Theorem, $\left\{ {\frac{\mathcal{L}(z,t)}{a_{1}(t)}}\right\}
_{t\in I}$ is a normal family in $\mathcal{U}_{r_{3}}.$ From the analyticity
of $\frac{\partial \mathcal{L}(z,t)}{\partial t},$ we obtain that for all
fixed numbers $T>0$ and $r_{4},\;0<r_{4}<r_{3},$ there exists a constant $%
K_{1}>0$ (that depends on $T$ and $r_{4})$ such that%
\begin{equation*}
\left\vert {\frac{\partial \mathcal{L}(z,t)}{\partial t}}\right\vert <K_{1},%
\text{ \ \ }\forall z\in \mathcal{U}_{r_{4}},\;t\in \left[ {0,T}\right] .
\end{equation*}%
Therefore, the function $\mathcal{L}(z,t)$ is locally absolutely continuous
in $I,$ locally uniformly with respect to $\mathcal{U}_{r_{4}}.$

Consider the function $p:\mathcal{U}_{r}\times I\rightarrow \mathbb{C}$ for$%
\;0<r<r_{4}$ and $t\in I,$ defined by%
\begin{equation*}
p(z,t)={z\frac{\partial \mathcal{L}(z,t)}{\partial z}}\diagup \frac{\partial 
\mathcal{L}(z,t)}{\partial t}.
\end{equation*}%
The function $p(z,t)$ is analytic in $\mathcal{U}_{r},$ $0<r<r_{4}.$ If the
function%
\begin{equation}
w(z,t)=\frac{p(z,t)-1}{p(z,t)+1}=\frac{\frac{z\partial \mathcal{L}(z,t)}{%
\partial z}-\frac{\partial \mathcal{L}(z,t)}{\partial t}}{\frac{z\partial 
\mathcal{L}(z,t)}{\partial z}+\frac{\partial \mathcal{L}(z,t)}{\partial t}}
\label{eq8}
\end{equation}%
\noindent is analytic in $\mathcal{U}\times I$ and $\left\vert {w(z,t)}%
\right\vert <1,$ for all $z\in \mathcal{U}\;$and $t\in I,$ then $p(z,t)$ has
an analytic extension with positive real part in $\mathcal{U},$ for all $%
t\in I.$ From equality (\ref{eq8}) we have%
\begin{equation}
w(z,t)=\frac{(1+a)\mathcal{G}(z,t)+1-b}{(1-a)\mathcal{G}(z,t)+1+b},
\label{eq9}
\end{equation}%
\noindent where%
\begin{eqnarray}
\mathcal{G}(z,t) &=&e^{-(a+b)t}\left\{ {\left( {\frac{(1+c){f}^{\prime }({%
e^{-at}}z)}{[\mathcal{R}^{n}h({e^{-at}}z){]}^{\prime }-\alpha }-1}\right) }%
\right.  \label{eq10} \\
&&\left. +{(e^{(a+b)t}-1)\left[ {(\beta -1)\frac{{e^{-at}}z{g}^{\prime }({%
e^{-at}}z)}{g({e^{-at}}z)}+\frac{{e^{-at}}z[\mathcal{R}^{n}h({e^{-at}}z){]}%
^{\prime \prime }}{[\mathcal{R}^{n}h({e^{-at}}z){]}^{\prime }-\alpha }}%
\right] }\right\}  \notag
\end{eqnarray}%
for $z\in \mathcal{U}$ and $t\in I.$

The inequality $\left\vert {w(z,t)}\right\vert <1$ for all $z\in \mathcal{U}%
\;$and $t\in I,$ where $w(z,t)$ is defined by (\ref{eq9}), is equivalent to%
\begin{equation}
\left\vert {\mathcal{G}(z,t)-\frac{m-1}{2}}\right\vert <\frac{m+1}{2},\text{
\ \ }\forall z\in \mathcal{U},\;t\in I.  \label{eq11}
\end{equation}%
Let us denote%
\begin{equation}
\mathcal{H}(z,t)=\mathcal{G}(z,t)-\frac{m-1}{2},\text{ \ \ }\forall z\in 
\mathcal{U},\;t\in I.  \label{eq111}
\end{equation}%
In view of (\ref{eq41}), (\ref{eq5}), from (\ref{eq10}) and (\ref{eq111}) we
have%
\begin{equation}
\left\vert {\mathcal{H}(z,0)}\right\vert =\left\vert {\left( {\frac{(1+c){f}%
^{\prime }(z)}{[\mathcal{R}^{n}h(z){]}^{\prime }-\alpha }-1}\right) -\frac{%
m-1}{2}}\right\vert <\frac{m+1}{2},  \label{eq12}
\end{equation}%
and for $z=0,$ $t>0,$ since (\ref{eq411})%
\begin{eqnarray*}
\left\vert {\mathcal{H}(0,t)}\right\vert &=&\left\vert {e^{-(a+b)t}}\left( {%
\frac{c+\alpha }{1-\alpha }}\right) {+(1-e^{-(a+b)t})(\beta -1)-\frac{m-1}{2}%
}\right\vert \\
&=&\left\vert {e^{-(a+b)t}}\left( {\frac{1+c}{1-\alpha }-\frac{m+1}{2}}%
\right) {+(1-e^{-(a+b)t})}\left( {\beta -\frac{m+1}{2}}\right) \right\vert \\
&<&{e^{-(a+b)t}}\frac{m+1}{2}{+(1-e^{-(a+b)t})}\frac{m+1}{2}=\frac{m+1}{2}.
\end{eqnarray*}%
Since $\left\vert {e^{-at}z}\right\vert \leqslant \left\vert {e^{-at}}%
\right\vert =e^{-at}<1$ for all $z\in \overline{\mathcal{U}}=\left\{ {z\in 
\mathbb{C}:\;\left\vert z\right\vert \leqslant 1}\right\} $ and $t>0,$ we
find that $\mathcal{H}(z,t)$ is an analytic function in $\overline{\mathcal{U%
}}.$ Making use of the maximum modulus principle we obtain that for each $%
t>0 $ arbitrarily fixed there exists $\theta =\theta (t)\in \mathbb{R}$ such
that%
\begin{equation}
\left\vert {\mathcal{H}(z,t)}\right\vert <\underset{\left\vert \zeta
\right\vert =1}{\max }\left\vert {\mathcal{H}(\zeta ,t)}\right\vert
=\left\vert {\mathcal{H}(e^{i\theta },t)}\right\vert ,  \label{eq13}
\end{equation}%
\noindent for all $z\in \mathcal{U}.$

Let us denote $u=e^{-at}e^{i\theta }.$ Then $\left\vert u\right\vert
=e^{-at},$ $e^{-(a+b)t}=(e^{-at})^{m+1}=\left\vert u\right\vert ^{m+1}$ and
from (\ref{eq10}) we have%
\begin{eqnarray*}
\left\vert {\mathcal{H}(e^{i\theta },t)}\right\vert &=&\left\vert {%
\left\vert u\right\vert ^{m+1}\left( {\frac{(1+c){f}^{\prime }(u)}{[\mathcal{%
R}^{n}h(u){]}^{\prime }-\alpha }-1}\right) }\right. \\
&&+\left. {(1-\left\vert u\right\vert ^{m+1})\left[ {(\beta -1)\frac{u{g}%
^{\prime }(u)}{g(u)}+\frac{u[\mathcal{R}^{n}h(u){]}^{\prime \prime }}{[%
\mathcal{R}^{n}h(u){]}^{\prime }-\alpha }}\right] -\frac{m-1}{2}}\right\vert
.
\end{eqnarray*}%
Since $u\in U,$ the inequality (\ref{eq5}) implies that%
\begin{equation}
\left\vert {\mathcal{H}(e^{i\theta },t)}\right\vert \leqslant \frac{m+1}{2},
\label{eq14}
\end{equation}%
\noindent and from (\ref{eq12}) and (\ref{eq14}) it follows that the
inequality (\ref{eq11})%
\begin{equation*}
\left\vert {\mathcal{H}(z,t)}\right\vert =\left\vert {\mathcal{G}(z,t)-\frac{%
m-1}{2}}\right\vert <\frac{m+1}{2}
\end{equation*}%
\noindent is satisfied for all $z\in \mathcal{U}\;$and $t\in I.$ Therefore $%
\left\vert {w(z,t)}\right\vert <1$ for all $z\in \mathcal{U}\;$and $t\in I.$

Since all the conditions of Theorem \ref{T1} are satisfied, we obtain that
the function $\mathcal{L}(z,t)$ has an analytic and univalent extension to
the whole unit disk $\mathcal{U},$ for all $t\in I.$ For $t=0$ we have $%
\mathcal{L}(z,0)=\mathcal{F}_{\beta }(z),$ for $z\in \mathcal{U}$ and
therefore, the function $\mathcal{F}_{\beta }(z)$ is analytic and univalent
in $\mathcal{U}.$
\end{proof}

For $n=0,$ in Theorem \ref{T2} we obtain another univalence criterion as
follows.

\begin{corollary}
\label{C2}Let $f,g,h\in \mathcal{A}.$ Also \textit{let} $m\in \mathbb{R}_{+}$%
\textit{\ and complex numbers }$\alpha ,$ $\beta ,$ $c$\textit{\ be }as in
Theorem \ref{T2}. \textit{If the inequalities }%
\begin{equation}
\left\vert {\left( {\frac{(1+c){f}^{\prime }(z)}{{h}^{\prime }(z)-\alpha }-1}%
\right) -\frac{m-1}{2}}\right\vert <\frac{m+1}{2}  \label{eq141}
\end{equation}%
\textit{and}%
\begin{equation*}
\left\vert {\left\vert z\right\vert ^{m+1}\left( {\frac{(1+c){f}^{\prime }(z)%
}{{h}^{\prime }(z)-\alpha }-1}\right) }\right.
\end{equation*}%
\begin{equation}
\left. {+(1-\left\vert z\right\vert ^{m+1})\left[ {(\beta -1)\frac{z{g}%
^{\prime }(z)}{g(z)}+\frac{z{h}^{\prime \prime }(z)}{{h}^{\prime }(z)-\alpha 
}}\right] -\frac{m-1}{2}}\right\vert \leqslant \frac{m+1}{2}  \label{eq15}
\end{equation}%
\textit{are true for all }$z\in \mathcal{U},$\textit{\ then the function }$%
\mathcal{F}_{\beta }(z)$ defined by (\ref{eq6}) \textit{is analytic and
univalent in }$\mathcal{U},$\textit{\ where the principal branch is intended.%
}
\end{corollary}

For $h(z)=f(z)$ and $\alpha =m-1=0$ in Corollary \ref{C2} we obtain
following criterion which will be useful for univalency of the general
integral operators.

\begin{corollary}
\label{C1} \textit{Let }$c,\beta $\textit{\ be complex numbers such that }$%
\left\vert c\right\vert <1,$ $\left\vert {\beta -1}\right\vert <1$\textit{\
and let the functions }$f,g\in \mathcal{A}.$\textit{\ If the inequality}%
\begin{equation*}
\left\vert {c\left\vert z\right\vert ^{2}+(1-\left\vert z\right\vert ^{2})%
\left[ {(\beta -1)\frac{z{g}^{\prime }(z)}{g(z)}+\frac{z{f}^{\prime \prime
}(z)}{{f}^{\prime }(z)}}\right] }\right\vert \leqslant 1
\end{equation*}%
\textit{is true for all }$z\in \mathcal{U},$\textit{\ then the function }$%
\mathcal{F}_{\beta }(z)$ defined by (\ref{eq6}) \textit{is analytic and
univalent in }$\mathcal{U},$\textit{\ where the principal branch is intended.%
}

\begin{remark}
\label{R2} For special values of parameters $\beta ,c$ and $m$ our results
reduce to several well-known results as follows:

\begin{enumerate}
\item Putting $\beta =1$ and $\alpha =0$ in Corollary \ref{C1}, then we
obtain the the result of Ahlfors \cite{Ah}. For $c=0,$ Ahlfors's criterion
reduces to a criterion found earlier by Becker \cite{Be}.

\item Putting $h(z)=f(z),$ $\beta =m=1$ and $c=0$ in Corollary \ref{C2}, we
obtain the result of Pascu \cite{Pa}.

\item Putting $h(z)=f(z),$ $\beta =1$ and $c=0$ in Corollary \ref{C2}, we
obtain the result of R\u{a}ducanu \textit{et al}. \cite{RaOrDe}
\end{enumerate}
\end{remark}
\end{corollary}

For $h(z)=f(z)$ and $n=1,$ in Theorem \ref{T2}, we have the following
corollary.

\begin{corollary}
\label{C3}\textbf{\ }Let $f,g\in \mathcal{A}.$ Also \textit{let} $m\in 
\mathbb{R}_{+}$\textit{\ and complex numbers }$\alpha ,$ $\beta ,$ $c$%
\textit{\ be }as in Theorem \ref{T2}. \textit{If the inequalities }%
\begin{equation}
\left\vert {\left( {\frac{2(1+c)}{m+1}-1}\right) -\frac{z{f}^{\prime \prime
}(z)-\alpha }{{f}^{\prime }(z)}}\right\vert <\left\vert {1-\frac{z{f}%
^{\prime \prime }(z)-\alpha }{{f}^{\prime }(z)}}\right\vert  \label{eq151}
\end{equation}%
\textit{and }%
\begin{equation*}
\left\vert {\left\vert z\right\vert ^{m+1}\left( {\frac{(1+c){f}^{\prime }(z)%
}{z{f}^{\prime \prime }(z)+{f}^{\prime }(z)-\alpha }-1}\right) }\right.
\end{equation*}%
\begin{equation*}
\left. {+(1-\left\vert z\right\vert ^{m+1})\left[ {(\beta -1)\frac{z{g}%
^{\prime }(z)}{g(z)}+\frac{z^{2}{f}^{\prime \prime \prime }(z)+2z{f}^{\prime
\prime }(z)}{z{f}^{\prime \prime }(z)+{f}^{\prime }(z)-\alpha }}\right] -%
\frac{m-1}{2}}\right\vert \leqslant \frac{m+1}{2}
\end{equation*}%
\textit{are satisfied for all }$z\in \mathcal{U},$\textit{\ then the
function }$\mathcal{F}_{\beta }(z)$ defined by (\ref{eq6}) \textit{is
analytic and univalent in }$\mathcal{U},$\textit{\ where the principal
branch is intended.}

\begin{corollary}
\label{C4} \textit{Let }$f,h\in \mathcal{A}.$\textit{\ If the inequality}%
\begin{equation}
(1-\left\vert z\right\vert ^{2})\left\vert {\frac{z{f}^{\prime }(z)}{f(z)}+%
\frac{z{h}^{\prime \prime }(z)}{{h}^{\prime }(z)}}\right\vert \leqslant 1,
\label{eq17}
\end{equation}%
\textit{hold true for all }$z\in \mathcal{U},$\textit{\ then the function }$%
f(z)$\textit{\ is univalent in }$\mathcal{U}.$
\end{corollary}
\end{corollary}

\begin{proof}
It results from Theorem \ref{T2} with $n=c=\alpha =0;\;\beta =2;\;m=1\;$and $%
g(z)=f(z).$
\end{proof}

The condition (\ref{eq5}) of Theorem \ref{T2} can be replaced with a simpler
one.

\begin{theorem}
\label{T3} Let $f,g\in \mathcal{A}.$ Also \textit{let} $m\in \mathbb{R}_{+}$%
\textit{\ and complex numbers }$\alpha ,$ $\beta ,$ $c$\textit{\ be }as in
Theorem \ref{T2}.\textit{\ If the inequalities }%
\begin{equation}
\left\vert {\left( {\frac{(1+c){f}^{\prime }(z)}{[\mathcal{R}^{n}h(z){]}%
^{\prime }-\alpha }-1}\right) -\frac{m-1}{2}}\right\vert <\frac{m+1}{2}
\label{eq171}
\end{equation}%
\textit{and }%
\begin{equation}
\left\vert {(\beta -1)\frac{z{g}^{\prime }(z)}{g(z)}+\frac{z[\mathcal{R}%
^{n}h(z){]}^{\prime \prime }}{[\mathcal{R}^{n}h(z){]}^{\prime }-\alpha }-%
\frac{m-1}{2}}\right\vert \leqslant \frac{m+1}{2},  \label{eq18}
\end{equation}%
\textit{are true for all }$z\in \mathcal{U},$\textit{\ then the function }$%
\mathcal{F}_{\beta }(z)$ defined by (\ref{eq6}) \textit{is analytic and
univalent in }$\mathcal{U},$\textit{\ where the principal branch is intended}
\end{theorem}

\begin{proof}
Making use of (\ref{eq41}) and (\ref{eq18}) we obtain%
\begin{equation*}
\left\vert {\left\vert z\right\vert ^{m+1}\left( {\frac{(1+c){f}^{\prime }(z)%
}{[R^{n}h(z){]}^{\prime }-\alpha }-1}\right) +(1-\left\vert z\right\vert
^{m+1})\left[ {(\beta -1)\frac{z{g}^{\prime }(z)}{g(z)}+\frac{z[R^{n}h(z){]}%
^{\prime \prime }}{[R^{n}h(z){]}^{\prime }-\alpha }}\right] -\frac{m-1}{2}}%
\right\vert
\end{equation*}%
\begin{equation*}
=\left\vert {\left\vert z\right\vert ^{m+1}\left( {\frac{(1+c){f}^{\prime
}(z)}{[R^{n}h(z){]}^{\prime }-\alpha }-1-\frac{m-1}{2}}\right) }\right.
\end{equation*}%
\begin{equation*}
\left. {+(1-\left\vert z\right\vert ^{m+1})\left[ {(\beta -1)\frac{z{g}%
^{\prime }(z)}{g(z)}+\frac{z[\mathcal{R}^{n}h(z){]}^{\prime \prime }}{[%
\mathcal{R}^{n}h(z){]}^{\prime }-\alpha }-\frac{m-1}{2}}\right] }\right\vert
\end{equation*}%
\begin{equation*}
<\left\vert z\right\vert ^{m+1}\frac{m+1}{2}+(1-\left\vert z\right\vert
^{m+1})\frac{m+1}{2}=\frac{m+1}{2}
\end{equation*}%
The conditions of Theorem \ref{T2} being satisfied it follows that the
function $\mathcal{F}_{\beta }(z)$defined by (\ref{eq6}) is univalent in $%
\mathcal{U}.$
\end{proof}

For $n=c=\alpha =0,\;\beta =2,\;m=3,$ $h(z)=f(z)\;$and $g(z)=z,$ Theorem \ref%
{T3} yields

\begin{remark}
\label{R3} Let $f\in \mathcal{A}.$ If the inequality%
\begin{equation*}
\left\vert {\frac{z{f}^{\prime \prime }(z)}{{f}^{\prime }(z)}}\right\vert
\leqslant 2,
\end{equation*}%
\noindent is true for all $z\in \mathcal{U},$ then the function%
\begin{equation*}
\mathcal{F}_{2}(z)=\left( {2\int\limits_{0}^{z}{u{f}^{\prime }(u)du}}\right)
^{1\diagup 2}
\end{equation*}%
\noindent is analytic and univalent in $\mathcal{U}.$
\end{remark}

We consider $\alpha ,\;\beta \;$and $c$ be real numbers, such that $\alpha
<0.$ For $h(z)=f(z)$ and $n=0,$ by elementary calculations we obtain that
inequality (\ref{eq41}) is equivalent to%
\begin{equation}
\Re {f}^{\prime }(z)>\frac{m-c}{\alpha (m+1)}\left\vert {{f}^{\prime }(z)}%
\right\vert ^{2},\text{ \ }m\geq c;\text{\ }z\in \mathcal{U}.  \label{eq19}
\end{equation}

\begin{corollary}
\label{C5}\textbf{\ }Let $f\in \mathcal{A}.$ Also \textit{let} $m\in \mathbb{%
R}_{+}$\textit{\ and }$\alpha ,$ $\beta ,$ $c$\textit{\ be real numbers such
that}%
\begin{equation*}
\alpha <0,\text{ \ }c\neq -1,\;c\leq m,\ -1<c\leq m-\alpha (m+1),\text{ \ \ }%
\left\vert {\beta -\frac{m+1}{2}}\right\vert <\frac{m+1}{2}.
\end{equation*}%
\textit{If the inequalities}%
\begin{equation}
\Re {f}^{\prime }(z)>\frac{m-c}{\alpha (m+1)}\left\vert {{f}^{\prime }(z)}%
\right\vert ^{2}  \label{eq20}
\end{equation}%
and%
\begin{equation}
\left\vert {(\beta -1)\frac{z{f}^{\prime }(z)}{f(z)}+\frac{z{f}^{\prime
\prime }(z)}{{f}^{\prime }(z)-\alpha }-\frac{m-1}{2}}\right\vert \leqslant 
\frac{m+1}{2}  \label{eq21}
\end{equation}%
\textit{are satisfied for all }$z\in \mathcal{U},$\textit{\ then the
function }$f(z)$ \textit{is univalent in }$\mathcal{U},$\textit{\ where the
principal branch is intended.}
\end{corollary}

\begin{proof}
If we take the inequality (\ref{eq19}) and $f\in \mathcal{A}$ instead of the
inequality (\ref{eq171}) and $g,h\in \mathcal{A},$ respectively, and $n=0$
in Theorem \ref{T3}, Corollary \ref{C5} can be show easily.
\end{proof}

\begin{remark}
\label{R4} Consider $\beta =1$ and $\alpha <0$ in Corollary \ref{C5} If in
the inequality (\ref{eq20}) we let $\alpha \rightarrow -\infty $ we obtain
that%
\begin{equation*}
\Re {f}^{\prime }(z)>0,\text{ \ \ }z\in \mathcal{U}.
\end{equation*}%
Since (\ref{eq21}) holds true for $\beta =1$ and $\alpha \rightarrow -\infty 
$ it follows from Corollary \ref{C5} that the function $f$ \ is univalent in 
$\mathcal{U}.$ Therefore, we can conclude that the univalence criterion due
to \textit{Alexander-Noshiro-Warshawski } \cite{Al}, \cite{No}, \cite{Wa} is
a limit case of Corollary \ref{C5}.
\end{remark}

For suitable values of $\alpha ,\;\beta ,\;c,\;n,\;m,\;h(z)$ and $g(z)$ in
Theorem \ref{T3}, we can obtain different univelence criteria.

Reasoning along the same line as in the proof of the Theorem \ref{T1} for
the Loewner chain%
\begin{eqnarray}
&&\mathcal{L}(z,t)  \label{eq22} \\
&=&\left\{ {\beta \int\limits_{0}^{e^{-at}z}{g^{\beta -1}(s){f}^{\prime
}(s)ds+}}\frac{\beta {(e^{bt}-e^{-at})}}{1+c}{zg^{\beta -1}(e^{-at}z)\left( {%
\frac{\mathcal{R}^{n}f(e^{-at}z)-\alpha }{\mathcal{R}^{v}h(e^{-at}z)-\alpha }%
}\right) }\right\} ^{1\diagup \beta }  \notag
\end{eqnarray}%
we obtain the following theorem. We omit the details.

\begin{theorem}
\textit{\label{T4} }Let $f,g,h\in \mathcal{A}.$ Also \textit{let} $m\in 
\mathbb{R}_{+}$\textit{\ and }$\alpha ,$ $\beta ,$ $c$\textit{\ be complex
numbers such that}%
\begin{equation*}
\text{\ }c\neq -1,\;\ \left\vert c-\frac{m-1}{2}\right\vert \leqslant \frac{%
m+1}{2},\text{ \ \ }\left\vert {\beta -\frac{m+1}{2}}\right\vert <\frac{m+1}{%
2}.
\end{equation*}%
\textit{If the inequalities}%
\begin{equation}
\left\vert {\left( {(1+c){f}^{\prime }(z)\frac{\mathcal{R}^{v}h(z)-\alpha }{%
\mathcal{R}^{n}f(z)-\alpha }-1}\right) -\frac{m-1}{2}}\right\vert <\frac{m+1%
}{2}  \label{eq221}
\end{equation}%
\textit{and }%
\begin{equation*}
\left\vert {\left\vert z\right\vert ^{m+1}\left( {(1+c){f}^{\prime }(z)\frac{%
\mathcal{R}^{v}h(z)-\alpha }{\mathcal{R}^{n}f(z)-\alpha }-1}\right) }\right.
\end{equation*}%
\begin{equation}
+\left. {(1-\left\vert z\right\vert ^{m+1})\left[ {(\beta -1)\frac{z{g}%
^{\prime }(z)}{g(z)}+\frac{z[\mathcal{R}^{n}f(z){]}^{\prime }}{\mathcal{R}%
^{n}f(z)-\alpha }-\frac{z[\mathcal{R}^{v}h(z){]}^{\prime }}{\mathcal{R}%
^{v}h(z)-\alpha }}\right] -\frac{m-1}{2}}\right\vert \leqslant \frac{m+1}{2}
\label{eq23}
\end{equation}%
are \textit{true for all }$z\in \mathcal{U},$\textit{\ then the function }$%
\mathcal{F}_{\beta }(z)$ defined by (\ref{eq6}) \textit{is analytic and
univalent in }$\mathcal{U},$\textit{\ where the principal branch is intended.%
}
\end{theorem}

For $\alpha =0$ in Theorem \ref{T4}, we obtain new result as follows:

\begin{corollary}
\textit{\label{C51} }Let $f,g,h\in \mathcal{A}.$ Also \textit{let} $m\in 
\mathbb{R}_{+}$\textit{\ and complex numbers} $\beta ,$ $c$\textit{\ be }as
in Theorem \ref{T4}. \textit{If the inequalities}%
\begin{equation*}
\left\vert {\left( {(1+c){f}^{\prime }(z)\frac{\mathcal{R}^{v}h(z)}{\mathcal{%
R}^{n}f(z)}-1}\right) -\frac{m-1}{2}}\right\vert <\frac{m+1}{2}
\end{equation*}%
\textit{and }%
\begin{equation*}
\left\vert {\left\vert z\right\vert ^{m+1}\left( {(1+c){f}^{\prime }(z)\frac{%
\mathcal{R}^{v}h(z)}{\mathcal{R}^{n}f(z)}-1}\right) }\right.
\end{equation*}%
\begin{equation*}
+{(1-\left\vert z\right\vert ^{m+1})\left[ {(\beta -1)\frac{z{g}^{\prime }(z)%
}{g(z)}+(n+1)\frac{\mathcal{R}^{n+1}f(z)}{\mathcal{R}^{n}f(z)}-(v+1)\frac{%
\mathcal{R}^{v+1}h(z)}{\mathcal{R}^{v}h(z)}-n+v}\right] }
\end{equation*}%
\begin{equation*}
\left. {-\frac{m-1}{2}}\right\vert \leqslant \frac{m+1}{2}
\end{equation*}%
are \textit{true for all }$z\in \mathcal{U},$\textit{\ then the function }$%
\mathcal{F}_{\beta }(z)$ defined by (\ref{eq6}) \textit{is analytic and
univalent in }$\mathcal{U},$\textit{\ where the principal branch is intended.%
}
\end{corollary}

For $v=2,\;n=1$ and $v=0,\;n=2$ in Corollary \ref{C51} , we obtain Corollary %
\ref{C6} and Corollary \ref{C7}, respectively.

\begin{corollary}
\label{C6}\textbf{\ }Let $f,g,h\in \mathcal{A}.$ Also \textit{let} $m\in 
\mathbb{R}_{+}$\textit{\ and complex numbers }$\beta ,$ $c$\textit{\ be }as
in Theorem \ref{T4}. \textit{If the inequalities}%
\begin{equation}
\left\vert {(1+c)}\left( {{{2{h}^{\prime }(z)+}z{h}^{\prime \prime }(z)}}%
\right) {-m-1}\right\vert <m+1  \label{eq241}
\end{equation}%
\textit{and }%
\begin{equation*}
\left\vert {\left\vert z\right\vert ^{m+1}\left( {(1+c)\left( {{{{h}^{\prime
}(z)+}}}\frac{{z}}{2}{{{h}^{\prime \prime }(z)}}\right) -1}\right) }\right.
\end{equation*}%
\begin{equation}
+\left. {(1-\left\vert z\right\vert ^{m+1})\left[ {(\beta -1)\frac{z{g}%
^{\prime }(z)}{g(z)}+}\frac{zf^{\prime \prime }(z)}{f^{\prime }(z)}{-}\frac{%
3z^{2}h^{\prime \prime }(z)+z^{3}h^{\prime \prime \prime }(z)}{2zh^{\prime
}(z)+z^{2}h^{\prime \prime }(z)}\right] -\frac{m-1}{2}}\right\vert \leqslant 
\frac{m+1}{2},  \label{eq25}
\end{equation}%
\textit{are true for all }$z\in \mathcal{U},$\textit{\ then the function }$%
\mathcal{F}_{\beta }(z)$ defined by (\ref{eq6}) \textit{is analytic and
univalent in }$\mathcal{U},$\textit{\ where the principal branch is intended.%
}
\end{corollary}

\begin{corollary}
\label{C7}\textbf{\ }Let $f,g,h\in \mathcal{A}.$ Also \textit{let} $m\in 
\mathbb{R}_{+}$\textit{\ and complex numbers }$\beta ,$ $c$\textit{\ be }as
in Theorem \ref{T4}. \textit{If the inequalities}%
\begin{equation}
\left\vert {\left( 2{(1+c){f}^{\prime }(z)\frac{h(z)}{2z{f}^{\prime
}(z)+z^{2}{{{f}^{\prime \prime }(z)}}}-1}\right) -\frac{m-1}{2}}\right\vert <%
\frac{m+1}{2}  \label{eq251}
\end{equation}%
\textit{and }%
\begin{equation*}
\left\vert {\left\vert z\right\vert ^{m+1}\left( 2{(1+c){f}^{\prime }(z)%
\frac{h(z)}{2z{f}^{\prime }(z)+z^{2}{{{f}^{\prime \prime }(z)}}}-1}\right) }%
\right.
\end{equation*}%
\begin{equation}
+\left. {(1-\left\vert z\right\vert ^{m+1})\left[ 1+{(\beta -1)\frac{z{g}%
^{\prime }(z)}{g(z)}+\frac{3z^{2}f^{\prime \prime }(z)+z^{3}f^{\prime \prime
\prime }(z)}{2zf^{\prime }(z)+z^{2}f^{\prime \prime }(z)}-\frac{z{h}^{\prime
}(z)}{h(z)}}\right] -\frac{m-1}{2}}\right\vert \leqslant \frac{m+1}{2}
\label{eq26}
\end{equation}%
\textit{are true for all }$z\in \mathcal{U},$\textit{\ then the function }$%
\mathcal{F}_{\beta }(z)$ defined by (\ref{eq6}) \textit{is analytic and
univalent in }$\mathcal{U},$\textit{\ where the principal branch is intended.%
}
\end{corollary}

\begin{remark}
\label{R6}\textbf{\ }Putting $\beta =m=n=1$, $\alpha =c=v=0$ and $h(z)=f(z)$
in Theorem \ref{T4}, we get the result of Kanas and Lecko \cite{KaLe}.
\end{remark}

\begin{corollary}
\label{C8}\textbf{\ }Let $f,g\in \mathcal{A}.$ Also \textit{let} $m\in 
\mathbb{R}_{+}$\textit{\ and complex numbers }$\beta ,$ $c$\textit{\ be }as
in Theorem \ref{T4}. \textit{If the inequalities }%
\begin{equation}
\left\vert {{(1+c){f}^{\prime }(z)}-\frac{m+1}{2}}\right\vert <\frac{m+1}{2}
\label{eq261}
\end{equation}%
\textit{and }%
\begin{equation}
\left\vert {\left\vert z\right\vert ^{m+1}\left( {(1+c){f}^{\prime }(z)-1}%
\right) +(1-\left\vert z\right\vert ^{m+1})(\beta -1){\frac{z{g}^{\prime }(z)%
}{{g}(z)}}-\frac{m-1}{2}}\right\vert \leqslant \frac{m+1}{2},  \label{eq27}
\end{equation}%
\textit{are true for all }$z\in \mathcal{U},$\textit{\ then the function }$%
\mathcal{F}_{\beta }(z)$ defined by (\ref{eq6}) \textit{is analytic and
univalent in }$\mathcal{U},$\textit{\ where the principal branch is intended.%
}
\end{corollary}

\begin{proof}
It results from Theorem \ref{T4} with $\alpha \rightarrow \infty $ $.$
\end{proof}

\begin{corollary}
\bigskip \label{C9}\textbf{\ }Let $f,g\in \mathcal{A}.$ Also \textit{let} $%
m\in \mathbb{R}_{+}$\textit{\ and complex numbers} $\beta ,$ $c$\textit{\ be 
}as in Theorem \ref{T4}. \textit{If the inequalities }%
\begin{equation*}
\left\vert {{(1+c){f}^{\prime }(z)}-\frac{m+1}{2}}\right\vert <\frac{m+1}{2}
\end{equation*}%
\textit{and }%
\begin{equation*}
\left\vert {(\beta -1)\frac{z{g}^{\prime }(z)}{{g}(z)}-\frac{m-1}{2}}%
\right\vert \leqslant \frac{m+1}{2},
\end{equation*}%
\textit{are true for all }$z\in \mathcal{U},$\textit{\ then the function }$%
\mathcal{F}_{\beta }(z)$ defined by (\ref{eq6}) \textit{is analytic and
univalent in }$\mathcal{U},$\textit{\ where the principal branch is intended.%
}
\end{corollary}

\begin{proof}
Corollary \ref{C9} can be demonstrated the same line as in the proof of the
Theorem \ref{T3}. We omit the details.
\end{proof}

\section{\textbf{Univalence criteria and the S\u{a}l\u{a}gean derivative}}

S\u{a}l\u{a}gean (see \cite{Sal}) introduced an operator $S^{n}$ $(n\in 
\mathbb{N}_{0})$ defined, for a function $f\in \mathcal{A},$ by 
\begin{eqnarray*}
S^{0}f(z) &=&f(z) \\
S^{1}f(z) &=&Sf(z)=zf^{\prime }(z) \\
S^{2}f(z) &=&S(Sf(z))=zf^{\prime }(z)+z^{2}f^{\prime \prime }(z) \\
S^{n}f(z) &=&S(S^{n-1}f(z)).
\end{eqnarray*}

Replacing the Ruscheweyh derivative by the S\u{a}l\u{a}gean derivative in
the construction of the Loewner chain, applying Theorem \ref{T1} and using
the well known condition $z\left( S^{n}f(z)\right) ^{\prime }=S^{n+1}f(z)$
we obtain Theorems \ref{T5} and \ref{T6}.

\begin{theorem}
\label{T5} Let $f,g,h\in \mathcal{A}.$ Also \textit{let} $m\in \mathbb{R}%
_{+} $\textit{\ and }$\alpha ,$ $\beta ,$ $c$\textit{\ be complex numbers
such that}%
\begin{equation*}
\mathit{\ }\alpha \neq 1,\text{ \ }c\neq -1,\;\ \left\vert \frac{1+c}{%
1-\alpha }-\frac{m+1}{2}\right\vert \leqslant \frac{m+1}{2},\text{ \ \ }%
\left\vert {\beta -\frac{m+1}{2}}\right\vert <\frac{m+1}{2}.
\end{equation*}%
\textit{If the inequalities }%
\begin{equation*}
\left\vert {\left( {\frac{(1+c){f}^{\prime }(z)}{[S^{n}h(z){]}^{\prime
}-\alpha }-1}\right) -\frac{m-1}{2}}\right\vert <\frac{m+1}{2}
\end{equation*}%
\textit{and }%
\begin{equation*}
\left\vert {\left\vert z\right\vert ^{m+1}\left( {\frac{(1+c){f}^{\prime }(z)%
}{[S^{n}h(z){]}^{\prime }-\alpha }-1}\right) }\right.
\end{equation*}%
\begin{equation*}
\left. {+(1-\left\vert z\right\vert ^{m+1})\left[ {(\beta -1)\frac{z{g}%
^{\prime }(z)}{g(z)}+\frac{z[S^{n}h(z){]}^{\prime \prime }}{[S^{n}h(z){]}%
^{\prime }-\alpha }}\right] -\frac{m-1}{2}}\right\vert \leqslant \frac{m+1}{2%
}
\end{equation*}%
\textit{are true for all }$z\in \mathcal{U},$\textit{\ then the function }$%
\mathcal{F}_{\beta }(z)$ defined by (\ref{eq6}) \textit{is analytic and
univalent in }$\mathcal{U},$\textit{\ where the principal branch is intended.%
}
\end{theorem}

\begin{theorem}
\textit{\label{T6} }Let $f,g,h\in \mathcal{A}.$ Also \textit{let} $m\in 
\mathbb{R}_{+}$\textit{\ and }$\alpha ,$ $\beta ,$ $c$\textit{\ be complex
numbers such that}%
\begin{equation*}
c\neq -1,\;\ \left\vert c-\frac{m-1}{2}\right\vert \leqslant \frac{m+1}{2},%
\text{ \ \ }\left\vert {\beta -\frac{m+1}{2}}\right\vert <\frac{m+1}{2}.
\end{equation*}%
\textit{If the inequalities}%
\begin{equation*}
\left\vert {\left( {(1+c){f}^{\prime }(z)\frac{S^{v}h(z)-\alpha }{%
S^{n}f(z)-\alpha }-1}\right) -\frac{m-1}{2}}\right\vert <\frac{m+1}{2}
\end{equation*}%
\textit{and }%
\begin{equation*}
\left\vert {\left\vert z\right\vert ^{m+1}\left( {(1+c){f}^{\prime }(z)\frac{%
S^{v}h(z)-\alpha }{S^{n}f(z)-\alpha }-1}\right) }\right.
\end{equation*}%
\begin{equation*}
+\left. {(1-\left\vert z\right\vert ^{m+1})\left[ {(\beta -1)\frac{z{g}%
^{\prime }(z)}{g(z)}+\frac{S^{n+1}f(z)}{S^{n}f(z)-\alpha }-\frac{S^{v+1}h(z)%
}{S^{v}h(z)-\alpha }}\right] -\frac{m-1}{2}}\right\vert \leqslant \frac{m+1}{%
2},
\end{equation*}%
are \textit{true for all }$z\in \mathcal{U},$\textit{\ then the function }$%
\mathcal{F}_{\beta }(z)$ defined by (\ref{eq6}) \textit{is analytic and
univalent in }$\mathcal{U},$\textit{\ where the principal branch is intended.%
}
\end{theorem}

Different corollaries from Theorem \ref{T5} and \ref{T6} can be obtained
under suitable choices of $m,\beta ,v,\alpha ,n$ and $h(z),g(z).$ For
example, if $v=2,$ $n=1$ and $\alpha =0$ then Theorem \ref{T6} reduces to
Corollary \ref{C10}.

\begin{corollary}
\label{C10} Let $f,g\in \mathcal{A}.$ Also \textit{let} $m\in \mathbb{R}_{+} 
$\textit{\ and complex numbers} $\beta ,$ $c$\textit{\ be }as in Theorem \ref%
{T6}. \textit{If the inequalities}%
\begin{equation*}
\left\vert {{(1+c)}({h^{\prime }(z)+zh^{\prime \prime }(z))}-\frac{m+1}{2}}%
\right\vert <\frac{m+1}{2}
\end{equation*}%
\textit{and }%
\begin{equation*}
\left\vert {\left\vert z\right\vert ^{m+1}\left( {(1+c)({h^{\prime
}(z)+zh^{\prime \prime }(z))}-1}\right) }\right.
\end{equation*}%
\begin{equation*}
+\left. {(1-\left\vert z\right\vert ^{m+1})\left[ {(\beta -1)\frac{z{g}%
^{\prime }(z)}{g(z)}+}\frac{zf^{\prime \prime }(z)}{f^{\prime }(z)}{-}\frac{%
2z^{2}h^{\prime \prime }(z)+z^{3}h^{\prime \prime \prime }(z)}{zh^{\prime
}(z)+z^{2}h^{\prime \prime }(z)}\right] -\frac{m-1}{2}}\right\vert \leqslant 
\frac{m+1}{2},
\end{equation*}%
are \textit{true for all }$z\in \mathcal{U},$\textit{\ then the function }$%
\mathcal{F}_{\beta }(z)$ defined by (\ref{eq6}) \textit{is analytic and
univalent in }$\mathcal{U},$\textit{\ where the principal branch is intended.%
}
\end{corollary}

\begin{remark}
\label{R5} In its special case when $m,\;\beta ,\;v,\;\alpha ,\;n,$ $h(z)$
and $g(z)$ in Theorem \ref{T2} - \ref{T6} we obtain the some results of
Kanas \textit{and Srivastava }\cite{KaSr}.
\end{remark}

\textbf{Acknowledgement.} \textit{The present investigation was supported by
Atat\"{u}rk University Rectorship under BAP Project (The Scientific and
Research Project of Atat\"{u}rk University) Project No: 2010/28.}

\begin{center}
.
\end{center}

\end{document}